\documentclass[12pt]{elsarticle} 
\usepackage{graphicx}
\usepackage{amsmath}
\usepackage{bm}
\usepackage{amsfonts}
\usepackage{amsthm}

\newtheorem{Theorem}{Theorem}[section]
\newtheorem{Lemma}{Lemma}[section]
\newtheorem{Remark}{Remark}[section]

\numberwithin{equation}{section}
\def\ee{\end{eqnarray*}}
\def\be{\begin{eqnarray*}}
\def\bee{\end{eqnarray}}
\def\bbe{\begin{eqnarray}}
\def\ea{\end{align*}}
\def\ba{\begin{align*}}
\makeatletter
\let\today\relax
\def\ps@pprintTitle{%
    \let\@oddhead\@empty
    \let\@evenhead\@empty
    \def\@oddfoot{\footnotesize\itshape
      \hfill\today}
    \let\@evenfoot\@oddfoot
    }
\makeatother

\makeatletter
\newcommand{\bal}{\@ifstar{\@bals}{\@bal}}
\def\@bals#1\eal{\begin{align*}#1\end{align*}}
\def\@bal#1\eal{\begin{align}#1\end{align}}
\makeatletter

\def\u{\bm u}

\def\e{\bm e}

  \def\R{\mathbb{Re}}

    \def\e{\mathrm e}
    \def\T{ \mathbb{T}_\alpha}
    
    \def\p{\partial}

\begin{document}
\journal{}
\title{Enhanced  and unenhanced dampings of the Kolmogorov flow   }

\author{Zhi-Min Chen}

 \ead{zmchen@szu.edu.cn}
\address{School  of Mathematical Sciences, Shenzhen University, Shenzhen 518060,  China}%

\begin{abstract}

The Kolmogorov flow  represents  the stationary sinusoidal solution $(\sin y,0)$ to a two-dimensional spatially periodic Navier-Stokes system, driven by an external force.
This system admits the additional non-stationary solution  $(\sin y,0)+e^{-\nu t} (\sin y,0)$, which  tends exponentially to the Kolmogorov flow at the minimum decay rate  determined by the viscosity $\nu$. Enhanced damping or enhanced dissipation of the problem is obtained by presenting higher decay rate for the difference between a solution and the non-stationary basic solution. Moreover, for the understanding of the metastability problem in an explicit manner, a variety of exact solutions are presented to show enhanced and unenhanced dampings.

\end{abstract}
\begin{keyword}
Kolmogorov flow, two-dimensional Navier-Stokes equations, metastability, enhanced damping, enhanced dissipation

\MSC[2020] 35B35 \sep 35Q30 \sep 35Q35 \sep 76D05
\end{keyword}

\maketitle

\section{Introduction}

Consider a two-dimensional periodically driven flow with viscosity $\nu>0$  governed by the    Navier-Stokes equation system
\bbe \label{n1}\p_t \u + \u\cdot \nabla \u -\nu \Delta\u +\nabla p=\nu (a \sin y, 0),\,\,\,\nabla \cdot \u=0
\bee
in the  flat torus  $\T = [0, 2\pi/\alpha) \times [0,2\pi)$ for $\alpha>0$ and $a\ge 0$. The pressure  $p$  and velocity  $\u$ are   subject to the  spatially periodic  condition
with respect $2\pi/\alpha$ period in $x$ and $2\pi$ period in $y$.
The vorticity formulation  of (\ref{n1}) is expressed as
\be \p_t\omega + \u\cdot\nabla \omega -\nu \Delta\omega =- \nu a\cos y\ee
or
\bbe\label{NS} \p_t\omega +J(\Delta^{-1}\omega,\omega) - \nu \Delta\omega  =- \nu a\cos y\bee
with the Jacobian $J(\phi,\varphi) = \p_x \phi\p_y\varphi-\p_y\phi\p_x\varphi$. The velocity $\u$, vorticity $\omega$ and stream function $\psi$ of the flow  are subject to the identities
\be\omega= -\Delta \psi \,\,\mbox{ and } \,\,\, \u=(\p_y \psi, -\p_x\psi).
\ee

This problem was introduced by Kolmogorov  in 1959 in a seminar \cite{Arn}, where he presented the basic stationary flow
\bbe\label{Kol} \u= (\sin y,0) \mbox{ or } \omega = -\cos y,\bee
 which   is known as the Kolmogorov flow, and encouraged  turbulence study from the instability of the flow.

When $\alpha =1$, the linear stability of  the Kolmogorov flow  for $\nu>0$ was given by Mishalkin and Sinai \cite{MS} and its nonlinear global stability  was obtained by Marchioro \cite{March}.
When $0<\alpha < 1$, Iudovich \cite{Yud} showed the instability of the Kolmogorov flow leading to the occurrence of secondary  stationary  flows. For some $\alpha \in (0,1)$, the existence of Hopf bifurcation was obtained by Chen and Price \cite{Chen1999}. Moreover,  if $\frac{\sqrt{3}}4\le \alpha < \frac{\sqrt{3}}2$ and (\ref{NS}) is associated with  a  horizontal a free-slip boundary condition, perturbations around the Kolmogorov  flow give rise to   Hopf bifurcation (see  Chen and Price \cite{Chen2002}) into oscillatory flows. Nonlinear interaction arises from the coexistence of multiple  secondary oscillatory flows from Hopf bifurcation and  develops into chaotic flows (see Chen and Price \cite{Chen2005}).  Turbulent Kolmogorov flow was also studied by Chandler and  Kerswell \cite{Tur}.

In the present study, we are interested in  stability problem of the Kolmogorov flow. Especially, for the unforced case $a=0$,  the $L_2$ inner product of (\ref{NS}) with $\omega$  produces the  $L_2$ energy estimate
 \be \frac12\frac{d}{dt}\|\omega(t)\|_{L_2}^2 +\nu \|\nabla \omega (t)\|_{L_2}^2= -\int_{\T}J(\Delta^{-1} \omega,\omega)\omega dxdy =0.
 \ee
 This yields  the  damping  of the enstrophy
\bbe\label{decay} \|\omega(t)\|_{L_2} \le e^{-\nu t} \|\omega(0)\|_{L_2},\,\,\mbox{ for } \alpha >0, \bee
at the minimum decay rate $\nu$. However, this damping may  be  enhanced for some class of  initial data with respect to  small $\nu$ or large Reynolds number $(Re=1/\nu)$.
The enhanced damping problem was raised by   Beck and  Wayne \cite{BW}.
They considered the  equation
\bbe\label{n6} \p_t\omega =\nu \Delta\omega  - e^{-\nu t} \sin y \p_x (1+\Delta^{-1}) \omega , \bee
which is linearized from the Navier-Stokes equation (\ref{NS}) with $a=0$ around its  exact solution
\bbe \label{n4}\omega=- e^{-\nu t} \cos y.
\bee
 When  the non-local operator $\Delta^{-1}$ is not involved, they obtained  the existence of positive constants $C$ and $M$ satisfying the estimate
\bbe \label{n7}\|\omega \|_{L_2} \le C e^{-M \sqrt{\nu} t}\|\omega(0)\|_{L_2}, \,\,\, 0< t\le \frac{\tau}\nu
\bee
for  given $\tau>0$, provided $\nu$ is small.
  The estimate (\ref{n7}) was recently obtained by Wei and Zhang \cite{Zhang} and Wei {\it et al.} \cite{Zhang2} for the complete linear equation  (\ref{n6}).
When $a=1$, Ibrahim {\it et al.} \cite{Ib} considered  the equation
 \bbe\label{n8} \p_t\omega =\nu \Delta\omega  -  \sin y \p_x (1+\Delta^{-1}) \omega , \bee
linearization of (\ref{NS})  around (\ref{Kol}),  and obtained the damping estimate (\ref{n7}).

The Kolmogorov flow moves in a horizontal direction along its streamlines  $y= $ constants, and is called a parallel flow. For the stability and instability of a nonparallel  vortex flow moving in multiple directions, we refer to the work of  author \cite{Chen2019,Chen2021}.

 The present study is  motivated by the work of  Lin and Xu \cite{Lin}, where they noticed that  linearized the Navier-Stokes flow (\ref{n6}) is close to the linearized Euler flow
\bbe\label{LE1} \p_t\omega =- \sin y \p_x (1+\Delta^{-1}) \omega  \bee
for  $t$  moderate and $\nu$  sufficiently small. Thus  the low frequency modes leading to unenhanced damping in the sense of \cite{Lin} can be controlled by  the linearized Euler flow  due to a RAGE theorem \cite{Const} and the nonlinear convection term can be bounded  by small initial data.
They \cite{Lin} considered  the  enhanced damping property of the linear equation (\ref{n6}) and the nonlinear equation (\ref{NS}) with $a=0$ as
\be \Big\|P_{\ne0}\omega\Big(\frac\tau\nu \Big) \Big\|_{L_2} <\delta \|P_{\ne0}\omega (0)\|_{L_2} \ee
for any  $\tau,\,\delta>0$, provided $\nu$ is small. Here
 $P_{\ne0}$  is  the projection  operator
\bbe P_{\ne0}\omega =\omega -P_0\omega \,\,\mbox{ with } \,\, P_0 \omega =\frac \alpha{2\pi} \int^{2\pi/\alpha}_0\omega dx.\label{Pne0}\bee

Actually, equation (\ref{NS}) has the additional solution
\bbe\label{n9} \u = (a\sin y + e^{-\nu t} \sin y,\,0)\,\,\, \mbox{ or } \,\,\,\omega = -a\cos y -  e^{-\nu t}\cos y. \bee
 Linearizing   the Navier-Stokes equation (\ref{NS}) around (\ref{n9}), we have
\bbe \label{n10}\p_t\omega = \nu\Delta \omega - a\sin y \p_x (1+\Delta^{-1})\omega -e^{-\nu t}  \sin y \p_x (1+\Delta^{-1})\omega.
\bee

The present study is twofold. Firstly, we extend the analysis of Lin and Xu \cite{Lin} to  the forced Navier-Stokes equation (\ref{NS}). Secondly, to the understanding of the problem in a different direction, we  present exact solutions showing  enhanced and unenhanced damping properties different to that of \cite{Lin}  in an explicit manner.

An enhanced damping result  reads as follows.
\begin{Theorem}\label{th1} For  any $\alpha>1$,   $\tau>0$,  $\delta >0$ and  $a\in \{0, \,1\}$, then the following assertions hold true.

 (i) Any  solution $\omega$ to  the linear equation (\ref{n10}) initiated from  $\omega(0)  \in L_2(\T)$ with $P_{\ne0}\omega(0) =\omega(0)$
satisfies
\be \Big\|\omega\Big(\frac\tau\nu \Big) \Big\|_{L_2} \le\delta \|\omega (0)\|_{L_2},
\ee
provided that $\nu>0$ is sufficiently small.

 (ii)For any  solution $\omega$ to  the nonlinear equation  (\ref{NS})  presented in the perturbed form
 \bbe \label{per}\omega = -a\cos y -  e^{-\nu t}\cos y+\omega',\bee
 the perturbed flow $\omega'$
satisfies
\bbe \label{aa0}\Big\|P_{\ne0}\omega'(\frac\tau\nu ) \Big\|_{L_2} \le\delta \Big\|P_{\ne0}\omega'(0)\Big\|_{L_2},
\bee
provided that $\omega'(0)\in L_2(\T)$,
\bbe \|\omega'(0) \|_{L_2} \le \nu d \label{aa1}
\bee
 for  $\nu$   and $d$   sufficiently small.
\end{Theorem}

When $a=1$, this result is on the  stability of  Kolmogorov flow. Hence metastability  rather than enhanced damping may be  more suitable
to address the present study.

When $a=0$, this  is comparable with an enhanced damping result    given by Lin and Xu \cite{Lin} in the following.

\begin{Theorem} (\cite[Theorem 1.2]{Lin})\label{th2}
Consider  the nonlinear Navier-Stokes equation (1.1) ($a=0$) on $\T$ with $\alpha\ge 1$. Denote
$P_{\mathcal K}$ to be the projection of $L_2(\T)$ to the subspace of Kolmogorov flows $W_{\mathcal K}={\rm span}\{\cos y, \sin y\}.$ Then,

 (i) (Rectangular torus) Suppose $\alpha>1$. There exists $d > 0$, such that for any
$\tau>0$ and $\delta >0$, if $\nu$ is small enough, then for any solution $\omega(t)$  to (\ref{NS}) with
initial vorticity $\omega(0) \in L_2(\T)$ satisfying
\bbe \|(I-P_{\mathcal K})\omega (0) \|_{L_2} \le d\nu, \label{Lin0}\bee
we have
\bbe \label{Lin1}\Big\|P_{\ne0} \omega\Big(\frac\tau \nu\Big)\Big\|_{L_2} < \delta  \Big\|P_{\ne0} \omega(0)\Big\|_{L_2}.\bee

(ii) (Square torus) Suppose $\alpha=1$. There exist $d > 0$, such that: for any
$M > 0$, $\tau >0$ and $\delta >0$, if $\nu$ is small enough , then for any solution $\omega(t)$  to (\ref{NS}) with initial data $\omega(0) \in L_2(\T)$ satisfying (\ref{Lin0}), 
either
\bbe \label{Lin2}\max_{ 0\le t\le \frac \tau\nu} \|P_a\omega(t)\|_{L_2} \ge M\|P_{\ne0}\omega(0)\|_{L_2}\bee
or
\bbe\label{Lin3} \inf_{0\le  t\le \frac\tau\nu }  \|(I-P_a)P_{\ne0}\omega(t)\|_{L_2} <\delta \|P_{\ne0}\omega(0)\|_{L_2}\bee
must hold true. Here $P_a$ is the orthogonal projection mapping $L_2$ onto the space span$\{\sin x,\cos x\}$.
\end{Theorem}

\begin{Remark}
It should be noted that Theorem \ref{th2} or  \cite[Theorem 1.2]{Lin} contradicts to  the following fact:

For $\tau=M=1$, $\delta =e^{-\alpha^2-1}$ and any constant $d>0$, the  function
\bbe \label{f}\omega=  d\nu \frac{\sqrt{\alpha}}{\sqrt{2}\pi}e^{-(\alpha^2+1)\nu t} \sin(\alpha x+y),\,\,\alpha\ge 1,
\bee
solves the Navier-Stokes equation (\ref{NS}) ($a=0$) in $\T$ and has the following properties:
\be \|(I-P_{\mathcal K})\omega(0)\|_{L_2}&=&\|\omega(0)\|_{L_2}=d\nu,\,\, \alpha \ge 1,
\\
P_{\ne0}\omega &=&\omega,\,\,\alpha \ge 1,
\\
 P_a \omega &=&0, \,\,\alpha =1,
\\
\|P_{\ne0}\omega(t)\|_{L_2} &\ge& \delta \|P_{\ne0}\omega(0)\|_{L_2},\,\, t \le \frac{\tau}{\nu}, \alpha >1,
\\
\|(I-P_{a})P_{\ne0}\omega(t)\|_{L_2} &\ge &\delta \|P_{\ne0}\omega(0)\|_{L_2},\,\, t \le \frac{\tau}{\nu}, \alpha =1.
\ee
Therefore,  the initial condition (\ref{Lin0}) holds ture, but none of the  conclusion estimates (\ref{Lin1})-(\ref{Lin3}) are valid, no matter how small the viscosity  $\nu$ is.

\end{Remark}

For a solution $\omega$ to the Navier-Stokes equation (\ref{NS}) with $a=0$ in the perturbed form
\bbe  \omega(t)=-\e^{-\nu t}\cos y + \omega'(t),\bee
we see that $P_{\ne0}\e^{-\nu t}\cos y =0$ and $(I-P_{\mathcal K})\e^{-\nu t}\cos y =0$. Hence there holds the identity
\bbe P_{\ne0}\omega(t)= P_{\ne0}\omega'(t).
\bee
Thus one might consider   $\omega'$ and $\omega$ to be the same in the enhanced damping analysis. However it should be noted that  the linear problem (\ref{n10}) is invariant in the subspace $(I-P_{\mathcal K})L_2(\T)$, but  this invariance cannot be extended to the nonlinear problem (\ref{NS}). Therefore the smallness assumption on $(I-P_{\mathcal K})\omega(0)$ is not enough to produce the desired enhanced damping estimates such as  (\ref{Lin1}).

More specifically, for $\omega$ representing the exact solution  (\ref{f}),  we see that
\be \omega(t) = -\e^{-\nu t} \cos y+\omega'(t)\ee
with
\be \omega'(t) = \omega(t) + \e^{-\nu t} \cos y = d\nu \frac{\sqrt{\alpha}}{\sqrt{2}\pi}e^{-(\alpha^2+1)\nu t} \sin(\alpha x+y)+ \e^{-\nu t} \cos y.
\ee
and hence
\bbe \|P_{\ne0}\omega(t)\|_{L_2}= \|P_{\ne0}\omega'(t)\|_{L_2}= \e^{-(\alpha^2+1)\nu t}d\nu,
\bee
but the initial functions  $\omega'(0)$ and $(I-P_{\mathcal K})\omega (0)$ are quit different as shown in the following.
\bbe \|\omega'(0)\|_{L_2}^2 &=&\| d\nu \frac{\sqrt{\alpha}}{\sqrt{2}\pi} \sin(\alpha x+y)\|_{L_2}^2+  \|\cos y\|_{L_2}^2\nonumber
\\
&=&  d^2\nu^2 + \frac{2\pi^2}\alpha>\frac{2\pi^2}\alpha,
\bee
which is  larger than the lower bound $\frac{2\pi^2}\alpha$  for any  viscosity $\nu>0$. In contrast, we see that
\be \|(I-P_{\mathcal K})\omega (0) \|_{L_2}&=&\|(I-P_{\mathcal K})\omega' (0) \|_{L_2}\nonumber
\\
&=& \|d\nu \frac{\sqrt{\alpha}}{\sqrt{2}\pi} \sin(\alpha x+y) \|_{L_2}=\nu d.
\ee
Hence to exclude the exact solution (\ref{f}) and to ensure the validity of Theorem \ref{th2}, Assumption (\ref{Lin0}) should be replaced by  (\ref{aa1}), or
\bbe \| \omega(0) +\e^{-\nu t}\cos y \|_{L_2} \le \nu d.
\bee

The enhanced damping property given by Theorem \ref{th1}  is determined by the estimates (\ref{aa0}) and (\ref{aa1}), which are not satisfied by (\ref{f}).
For compensating this shortness,  a variety of  exact solutions of (\ref{NS})  will be provided for showing their explicit  enhanced damping different to that presented in Theorem \ref{th1}.

Theorem \ref{th1}  involves only  the rectangular case $\alpha >1$, as the limit case $\alpha =1$  can be treated in a similar manner.

This paper is organized as follows.  Theorem \ref{th1} is to be  proved in Section 2. Exact solutions will be discussed in  Section 3.

\section{Proof of Theorem \ref{th1}}
The proof is  developed from    Lin and Xu \cite{Lin} together with  Constantin {\it et al.}\cite{Const}. Let $C$ be a generic constant and the condition $\alpha >1$ is assumed throughout this section.

For convenience of notation, we always assume that  $L_2(\T)$ denotes the $L_2$ space of all mean zero  functions:
\be L_2(\T)= \Bigg\{ \omega \,\,\Big|\,\, \int_{\T}\omega dxdy =0, \| \omega\|_{L_2} ^2 =\int_{\T} |\omega|^2 dxdy <\infty\Bigg\}.\ee

Recalling the operator $P_{\ne0}$ from (\ref{Pne0}) , we define the $L_2$ subspaces  $X= P_{\ne 0}L_2(\T)$, that is,
\be X =\Big\{ \omega  \in L_2(\T)\,\, |\,\, P_{\ne0}\omega=\omega\Big\}.\ee
 Note  that $(1+\Delta ^{-1})W_{\mathcal K}=\{0\}$. It is convenient to use the equivalent $L_2$-norm
\bal \|\omega \|_X = \|(1+\Delta^{-1})^{1/2}\omega\|_{L_2}
\eal
in the $L_2$ subspace  $(1+\Delta^{-1})L_2(\T)$.

\subsection{Preliminary estimates}

As given by Lin and Xu \cite{Lin}, the
enhanced damping lies on a modified  RAGE theorem \cite{Const} with respect to the linearized Euler equation
\bbe \label{nn1}\label{LE} \p_t \omega = -(a+1)\sin y \p_x(1+\Delta^{-1})\omega \bee
in $X$, by avoiding  unenhanced  damping from  low frequency eigenfunctions of the Laplacian operator $-\Delta$.
We adopt the complete eigenvalues
$$\alpha^2 \le \lambda_1\le \lambda_2 \le ... $$
 and the corresponding eigenfunctions $e_1$, $e_2,...$ of the operator $-\Delta$ in $X$. Let  $P_N$ be the $L_2$  projection mapping $X$ onto the subspace
span$\{e_1,...,e_N\}$. 

\begin{Lemma}\label{L1}(\cite[Lemma 2.3]{Lin} and \cite[Lemma 3.2]{Const})
Let $K\subset S \equiv\{\varphi \in X |\, \|\varphi \|_{L_2}=1\}$ be a compact set of $X$.
For any $N, \kappa >0$, there exists $T_c(N,\kappa,K)$ such that for all $T\ge T_c$ and $\omega^0(0)/\|\omega^0(0)\|_{L_2} \in K$,
\bbe \frac1T\int^T_0 \|P_N \omega^0(t) \|^2_{X}dt \le \kappa \|\omega^0(0)\|^2_{X},
\bee
where $\omega^0(t)$ is  the solution of the linearized Euler equation (\ref{nn1}) initiated from $\omega^0(0)$.
\end{Lemma}

Rewriting  the Navier-Stokes equation (\ref{NS})  perturbed around the basic flow   $-a\cos y -e^{-\nu t}\cos y$, we have
\bal\nonumber
\p_t  \omega^\nu -\nu \Delta \omega^\nu =&- a\sin y \p_x(1+\Delta^{-1}) \omega^\nu -  \e^{-\nu t } \sin y\p_x(1+\Delta^{-1}) \omega^\nu
\\
&-\sigma J(\Delta^{-1}\omega^\nu ,\omega^\nu ),\label{PNS}
\eal
with the constant $\sigma=0,\,1$. This equation is the linear  (\ref{n10})  for   $\sigma=0$ and is  the     nonlinear perturbed  form of (\ref{NS}) for $\sigma=1$.

The proof  of Theorem \ref{th1} is essentially based on  the fundamental $L_2$ energy  estimate  of   Kolmogorov problem (\ref{PNS}) given by Iudovich \cite[Eq. (2.14)]{Yud}
expressed as
\bal\label{L_2est} \frac{d}{dt}\int_{\T} \omega^\nu (t) (1+\Delta^{-1})\omega^\nu (t) dxdy  -2\nu \int_{\T} \Delta \omega^\nu (t)  (1+\Delta^{-1}) \omega^\nu (t) dxdy=0,
\eal
or
\bal\label{L_2est1} \frac{d}{dt} \|\omega^\nu (t) \|^2_X +2\nu\|\nabla \omega^\nu (t)\|^2_X=0.
\eal
This is imply due to the identities
  \be \int_{\T} (a+ \e^{-\nu t } )\sin y\p_x(1+\Delta^{-1}) \omega^\nu(1+\Delta^{-1})\omega^\nu  dxdy =0
  \ee
  and
\bbe  \int_{\T}  J(\Delta^{-1}\omega^\nu ,\omega^\nu )(1+\Delta^{-1})\omega^\nu  dxdy =0,\bee
after integration by parts.

When $a=1$, the energy estimate (\ref{L_2est}) gives rise to  the steady-state bifurcation analysis of \cite{Yud}  for $0<\alpha < 1$ and  the Hopf bifurcation analysis of \cite{Chen2002} for $\frac{\sqrt{3}}4\le \alpha <\frac{\sqrt{3}}2$.

As a consequence of  (\ref{L_2est1}), we have  an $L_2$ energy estimate of $P_{\ne0}\omega^\nu$  in the following.
\begin{Lemma}\label{newL1}Let $\omega^\nu$ be a  solution of (\ref{PNS}) with $\omega^\nu(0) \in X$. Assume that   $\| \omega^\nu(0)\|_{X}\le \nu d$ for a small constant $d>0$ when $\sigma=1$.  Then we have
\bal\label{nb1} \frac{d}{dt} \|P_{\ne 0}\omega^\nu(t) \|^2_X +\frac{15}8\nu\|\nabla P_{\ne 0}\omega^\nu(t)\|^2_X\le 0.
\eal
\end{Lemma}

When $a=0$, similar  estimate  has been shown in \cite[Lemma 3.2]{Lin}. For completion of analysis, we sketch a proof developed from \cite{Lin}.
\begin{proof}  Applying the operator $P_0$ to (\ref{PNS}), we have
\bal\nonumber
\p_t  P_0\omega^\nu(t)&-\nu \Delta P_0\omega^\nu(t)\\
=&- \sigma P_0 J(\Delta^{-1}P_{\ne0}\omega^\nu(t)+\Delta^{-1}P_{0}\omega^\nu(t),P_{\ne0}\omega^\nu(t)+P_{0}\omega^\nu(t))\nonumber\\
=&-\sigma P_0 J(\Delta^{-1}P_{\ne0}\omega^\nu(t), P_{\ne0}\omega^\nu(t)).\label{PNS0}
\eal
Taking the $L_2$ inner product of (\ref{PNS0}) with $2(1+\Delta^{-1})P_0\omega^\nu(t)$, we have
\bbe\lefteqn{ \frac{d}{dt} \| P_0\omega^\nu(t)\|^2_X+2\nu \|\nabla P_0 \omega^\nu(t)\|_X^2}\nonumber\\
 &=& -2\sigma \int_{\T} P_0 J(\Delta^{-1}P_{\ne0}\omega^\nu(t), P_{\ne0}\omega^\nu(t))(1+\Delta^{-1})P_0\omega^\nu(t) dxdy\nonumber
 \\
 &\ge & -2\sigma\|\nabla\Delta^{-1}P_{\ne0}\omega^\nu(t)\|_{L_\infty}\|\nabla P_{\ne0}\omega^\nu(t)\|_{L_2}\|(1+\Delta^{-1})P_0\omega^\nu(t)\|_{L_2}\nonumber
 \\
 &\ge & -2\sigma C\|\nabla P_{\ne0}\omega^\nu(t)\|_{X}^2\|\omega^\nu(t)\|_{X}\nonumber
 \\
 &\ge & -2\sigma C\|\nabla P_{\ne0}\omega^\nu(t)\|_{X}^2\|\omega^\nu(0)\|_{X},\label{new2}
\bee
after the use of the Sobolev imbedding theorem and (\ref{L_2est1}). Here and in what follows, $C$ is a generic constant independent of $\nu$.

On the other hand, rewriting (\ref{L_2est1}) as
\bal\nonumber \frac{d}{dt}& \!\|P_{\ne0}\omega^\nu(t) \|^2_X
\\ =&-2\nu\|\nabla P_{\ne0} \omega^\nu(t)\|^2_X  -\frac{d}{dt} \!\|P_{0}\omega^\nu(t) \|^2_X -2\nu\|\nabla P_{0} \omega^\nu(t)\|^2_X\label{L_2est11}
\eal
and taking (\ref{new2}) into account, we have
\bal
\frac{d}{dt} \!\|P_{\ne0}\omega^\nu(t) \|^2_X &\le -2(\nu-\sigma C\|\omega^\nu(0)\|_{X})\|\nabla P_{\ne0} \omega^\nu(t)\|^2_X\nonumber
\\
&\le -2\nu(1-\sigma C d)\|\nabla P_{\ne0} \omega^\nu(t)\|^2_X\le -\frac{15}8\nu \|\nabla P_{\ne0} \omega^\nu(t)\|^2_X,\nonumber
\eal
as $d$ is small so that  $2(1-Cd)>\frac{15}8$ when $\sigma =1$.
\end{proof}

The enhanced damping is to be confirmed  when a  solution to (\ref{PNS}) is close to  that of the linearized Euler equation (\ref{LE}) in the following sense.
\begin{Lemma}\label{L2} For  $\tau,\delta$, $d$ and $\nu$ given   in Theorem \ref{th1},
let $\omega^\nu$ be a  solution of (\ref{PNS}) with $\omega^\nu(0) \in L_2(\T)$  so that $\|\omega^\nu(0)\|_{L_2}<\nu d$ for $\sigma=1$ and let $\omega^0$ be a  solution of (\ref{LE}) with $\omega^0(0) \in X$ and $\|\nabla \omega^0(0)\|_X<\infty$.
If the estimate
\bbe \label{LL}\|P_{\ne0}\omega^\nu(t)\!-\!\omega^0(t)\|_X^2 \le \|P_{\ne0}\omega^\nu(0)\!-\!\omega^0(0)\|_X^2\!+\!C_1 \nu (1\!+\!t^3) \| \nabla\omega^0(0)\|_{X}^2\bee
 holds true  for a constant $C_1$ independent of $\nu$, 
 then
 \bbe \Big\|P_{\ne0}\omega^\nu\Big(\frac \tau\nu\Big) \Big\|_{L_2} < \delta \|P_{\ne0}\omega^\nu(0)\|_{L_2}.
\bee
\end{Lemma}

  For the unforced  case $a=0$, similar result   has been presented  in \cite[pp. 1824-1825]{Lin} originated from \cite[Proof of Theorem 1.4]{Const}.  Following  \cite{Const},
we  adopt   a   compact set $K\subset S$ required by  Lemma \ref{L1}.

\begin{proof}
For the  given $\delta$ and  $\tau>0$, we choose $N$ large enough  so that
\bbe e^{-\lambda_N \tau} <\frac{\alpha^2-1}{\alpha^2}\delta^2.\label{alpha}
\bee
In order to use Lemma \ref{L1}, we set  the compact set
\be K = \Big\{ \phi\in X |\,\, \|\nabla \phi\|^2_{X} \le \lambda_N, \,\, \|\phi\|_X=1 \Big\}.\ee
Let $t_1 = T_c(N,\frac1{10},K)$  from Lemma \ref{L1}  and $C_1$ from (\ref{LL}). Choose a small viscosity $\nu_0$ so that
\bbe\label{nn4} \nu_0 C_1(1+t_1^3) < \frac1{10\lambda_N} \,\,\,\mbox{ and } \,\,\,\nu_0 < \frac{\tau}{3t_1}.
\bee

For any $0<\nu<\nu_0$, if
\be \lambda_N \|P_{\ne0}\omega^\nu(t)\|^2_X\le \|\nabla P_{\ne0}\omega^\nu(t)\|^2_{X}
\ee
 for all $t\in [0,\frac \tau\nu ]$, it follows from  (\ref{nb1}) and (\ref{alpha})  that
\bbe
\Big\|P_{\ne0}\omega^\nu\Big(\frac\tau\nu\Big) \Big\|^2_X\le e^{-\lambda_N \tau} \|P_{\ne0}\omega^\nu(0)\|^2_X< \delta^2 \|P_{\ne0}\omega^\nu(0)\|^2_{L_2}\label{m1}
\bee
and we are done.

Otherwise,   there exits  $t_0\in [0, \frac \tau\nu )$  being the first time in the interval such that
\bbe \label{nn5}\|\nabla P_{\ne0}\omega^\nu(t_0)\|^2_{X} \le \lambda_N \|P_{\ne0}\omega^\nu(t_0)\|^2_X.
\bee

Now we take the solution  $\omega^0(t_0+t)$   of (\ref{LE})  for   $t\in  [0, t_1]$ initiated from
\bbe\label{vvv1}\omega^0(t_0)\equiv P_{\ne0}\omega^\nu(t_0)\in \{ \omega \in X|\, \|\nabla \omega\|_X <\infty\}.\bee
It follows from (\ref{LL}), (\ref{nn4}), (\ref{nn5}) and (\ref{vvv1}) that, for $t\in [0, t_1]$,
\bal \|P_{\ne0}\omega^\nu(t_0+t)-\omega^0(t_0+t)\|^2_X&\le \nu_0C_1(1+t_1^3) \|\nabla\omega^0(t_0)\|_{X}^2
\le \frac1{10} \|\omega^0(t_0)\|_X^2.\label{nn6}
\eal
Since  $\omega^0(t_0)/\|\omega^0(t_0)\|_X \in K$, by the definition of $t_1$ and Lemma \ref{L1}, we have
\bbe \label{nnb1}\frac1{t_1}\int^{t_1}_{0}\|P_N\omega^0(t_0+t) \|_X^2dt \le \frac1{10} \|\omega^0(t_0)\|_X^2.
\bee
By (\ref{nnb1}) and the conservation property  (\ref{L_2est1}) with $\nu=0$ for the linearized Euler equation (\ref{LE}),  we have
\bal* \frac1{t_1}\int^{t_1}_{0} \!\!\!\|(I-P_N)\omega^0(t_0+t) \|_X^2dt&=\|\omega^0(t_0+t)\|_X^2 -\frac1{t_1}\int^{t_1}_{0} \!\!\!\|P_N\omega^0(t_0+t) \|_X^2dt
\\
& =\|\omega^0(t_0)\|^2_X- \frac1{t_1}\int^{t_1}_{0}\|P_N\omega^0(t_0+t) \|_X^2dt
\\ &\ge  \frac9{10} \|\omega^0(t_0)\|_X^2.
\eal
This together with (\ref{nn6}) implies
\bal* \frac1{t_1}&\int^{t_1}_{0} \|(I-P_N)P_{\ne0}\omega^\nu(t_0+t) \|_X^2dt
\\
&\ge \frac1{t_1}\int^{t_1}_{0}\!\!\! \|(I-P_N)\omega^0(t_0+t) \|_X^2dt
-\frac1{t_1}\int^{t_1}_{0}\!\!\!\|P_{\ne0}\omega^\nu(t_0+t)-\omega^0(t_0+t) \|_X^2dt
\\
&\ge \frac45 \|\omega^0(t_0)\|_X^2.
\eal
Hence, we have
\bal*
\int^{t_1}_{0} \|\nabla P_{\ne0}\omega^\nu(t_0+t)\|^2_{X}dt &\ge \lambda_N\int^{t_1}_{0} \|(I-P_N)P_{\ne0}\omega^\nu(t_0+t)\|^2_X
\\
&\ge \frac{4\lambda_N t_1 }{5} \|\omega^0(t_0)\|_X^2.
\eal
By (\ref{nb1}) and the previous inequality, we have
\bal\|P_{\ne0}\omega^\nu(t_0+t_1)\|^2_X& \le  \|P_{\ne0}\omega^\nu(t_0)\|_{X}^2 -\frac{15}8\nu \int^{t_1}_{0}\|\nabla P_{\ne0}\omega^\nu(t_0+t)\|_{X}^2dt\nonumber
\\
& \le (1-\frac32\lambda_N \nu t_1 ) \|P_{\ne0}\omega^\nu(t_0)\|^2_{X} \nonumber
\\
&\le e^{-\frac32\lambda_N \nu t_1 } \|P_{\ne0}\omega^\nu(t_0)\|_X^2.\label{nn10}
\eal

Moreover,  for any interval $[\tau_0,\tau_1]$ such that
\bbe \|\nabla P_{\ne0}\omega^\nu(t)\|^2_X \ge \lambda_N \|P_{\ne0}\omega^\nu(t)\|^2_X,\,\,\,\, t\in [\tau_0,\tau_1],\bee
we have by (\ref{nb1}) that
\bbe \|P_{\ne0}\omega^\nu(\tau_1)\|^2_{X} \le e^{-\lambda_N \nu (\tau_1-\tau_0)} \|P_{\ne0}\omega^\nu(\tau_0)\|_X^2.\label{m1m}
\bee
Hence, the combination of all  decay estimates  from (\ref{nn10}) and  (\ref{m1m}), and  the use of $t_1<\frac{\tau}{3\nu}$ from (\ref{nn4}) imply the existence of $t_2 \in [\frac{2\tau}{3\nu}, \frac\tau\nu]$ such that
\bbe \label{vv1}\|P_{\ne0}\omega^\nu(t_2)\|^2_{X} \le e^{-\frac 32\lambda_N \nu  t_2} \|P_{\ne0}\omega^\nu(0)\|_X^2 \le e^{-\lambda_N  \tau} \|P_{\ne0}\omega^\nu(0)\|_X^2.
\bee
Therefore,  we have
\be
\Big\|P_{\ne0}\omega^\nu\Big(\frac\tau\nu\Big)\Big\|^2_{L_2}&\le& \frac{\alpha^2}{\alpha^2-1}\Big\|P_{\ne0}\omega^\nu\Big(\frac\tau\nu\Big)\Big\|^2_{X}
\\
& \le& \frac{\alpha^2}{\alpha^2-1}\|P_{\ne0}\omega^\nu(t_2)\|^2_{X}, \,\,\,\mbox{ by (\ref{nb1})},
\\
& \le& \frac{\alpha^2}{\alpha^2-1}e^{-\lambda_N \tau}\|P_{\ne0}\omega^\nu(0)\|^2_{X}, \,\,\,\mbox{ by (\ref{vv1})},
\\
&\le& \delta^2  \|P_{\ne0}\omega^\nu(0)\|^2_X, \,\,\,\mbox{ by (\ref{alpha})},
\\
& \le& \delta^2 \|P_{\ne0}\omega^\nu(0)\|_{L_2}^2.
\ee
This completes the proof of Lemma \ref{L2}.
\end{proof}

\subsection{Proof of Theorem \ref{th1}}

By Lemma \ref{L2}, it suffices to prove the validity of (\ref{LL}) for  the solutions $\omega^\nu$ of (\ref{PNS}) and $\omega^0$ of (\ref{LE}) when  $\omega^\nu(0) \in X$ and  $\omega^0(0) \in X^1$.

Let  $\omega=P_{\ne0}\omega^\nu-\omega^0$ . Thus the difference $\omega$ is subject to
the equation
\bal*
\p_t \omega
=& \nu \Delta P_{\ne0}\omega^\nu -  a\sin y \p_x(1+\Delta^{-1}) \omega-   \e^{-\nu t } \sin y\p_x(1+\Delta^{-1}) \omega
\\
&-(\e^{-\nu t }-1) \sin y\p_x(1+\Delta^{-1})\omega^0-\sigma P_{\ne0}J(\Delta^{-1}\omega^\nu,\omega^\nu).
\eal

Taking $L_2$ inner product of the previous equation with $(1+\Delta^{-1})\omega$ and integrating by parts, we have
\bal \nonumber
\frac12\frac{d}{dt}&\|\omega\|^2_X+\nu \| \nabla P_{\ne0}\omega^\nu\|_X^2
\\
=&  - \nu (\Delta P_{\ne0}\omega^\nu,  (1+\Delta^{-1})\omega^0)
 \nonumber
\\ \nonumber &-(\e^{-\nu t }-1) (\sin y\p_x(1+\Delta^{-1})\omega^0, (1+\Delta^{-1})P_{\ne0}\omega^\nu)
\\
&
-\sigma (P_{\ne0}J(\Delta^{-1}\omega^\nu, \omega^\nu),(1+\Delta^{-1}) \omega).\label{est}
\eal
According to integration by parts, the first two linear terms on the right-hand side of (\ref{est}) are bounded by
\bal
 \nu \|\nabla P_{\ne0}\omega^\nu\|_{X} & \|\nabla \omega^0\|_X +\nu t  \|\nabla (1+\Delta^{-1})\omega^0\|_{L_2}\|(1+\Delta^{-1}) P_{\ne0}\omega^\nu\|_{L_2}\nonumber
 \\
  &\le \nu \|\nabla P_{\ne0}\omega^\nu\|_{X}  \|\nabla \omega^0\|_X +\nu t  \| \omega^0\|_X\| P_{\ne0}\omega^\nu\|_X.\label{nnn2}
\eal
For the nonlinear problem $\sigma=1$, we adopt  the orthogonal decomposition
\be \omega^\nu=P_{\ne0}\omega^\nu+P_0\omega^\nu.
\ee
Integrating by parts and employing  H\"older inequality, Cauchy inequality and Sobolev imbedding, we  estimate  the third term on the right-hand side of (\ref{est}) as
\bal \nonumber
  -(P_{\ne0}J(&\Delta^{-1}\omega^\nu,\omega^\nu),(1+\Delta ^{-1})\omega)
  \\
=& -\Big(P_{\ne0}J(\Delta^{-1}P_0\omega^\nu,P_{\ne0}\omega^\nu)\!+\!P_{\ne0}J(\Delta^{-1}P_{\ne0}\omega^\nu,P_0\omega^\nu),(1\!+\!\Delta ^{-1})\omega\Big)
\nonumber
\\ &
-\Big(P_{\ne0}J(\Delta^{-1}P_{\ne0}\omega^\nu,P_{\ne0}\omega^\nu),(1+\Delta ^{-1})\omega\Big)
\nonumber
\\
 \le & \Big(\|\Delta ^{-1} P_0\omega^\nu\|_{L_\infty} \|\nabla P_{\ne0}\omega^\nu\|_{L_2} + \|\nabla \Delta^{-1}P_{\ne0}\omega^\nu \|_{L_\infty} \| P_0\omega^\nu\|_{L_2} \nonumber
 \\&+ \|  \Delta^{-1}P_{\ne0}\omega^\nu\|_{L_\infty} \|\nabla P_{\ne0}\omega^\nu \|_{L_2}\Big) \|\nabla  (1+\Delta ^{-1})\omega\|_{L_2}\nonumber
\\
\le & C\|\nabla P_{\ne0}\omega^\nu \|_{L_2} (\|P_0\omega^\nu\|_{L_2} + \|P_{\ne0}\omega^\nu \|_{L_2})\|\nabla (1+\Delta ^{-1})\omega\|_{L_2}\nonumber
\\
\le & C\|\nabla P_{\ne0}\omega^\nu \|_{X} \|\omega^\nu \|_{L_2}\|\nabla\omega\|_{X}\nonumber
\\
\le & C\|\nabla P_{\ne0}\omega^\nu \|_{X}^2 \|\omega^\nu \|_{L_2}+C\|\omega^\nu \|_{L_2}\|\nabla\omega^0\|_{X}^2
\label{nnn0}
\eal
Collecting terms of (\ref{est})-(\ref{nnn0}) and employing  Cauchy inequality, we have
\bal
\frac{d}{dt}\|\omega(t)\|^2_X\le& -2\nu \|\nabla P_{\ne0}\omega^\nu(t)\|_X^2
+2\nu \|\nabla P_{\ne0}\omega^\nu(t)\|_X  \|\nabla \omega^0(t)\|_X \nonumber
\\
&+2\nu t  \| P_{\ne0}\omega^\nu(t)\|_X\|\omega^0(t)\|_X
+C \sigma \|\nabla P_{\ne0}\omega^\nu(t)\|_{X}^2 \|\omega^\nu(t)\|_{L_2} \nonumber
\\
&+C\sigma\|\nabla P_{\ne0}\omega^\nu (t)\|_{X} \|\omega^\nu(t) \|_{L_2}\|\nabla\omega^0(t)\|_{X}\nonumber
\\
\le& -\nu \|\nabla P_{\ne0}\omega^\nu(t)\|_X^2
+2\nu\|\nabla \omega^0(t)\|_X^2 +2\nu t^2\| \omega^0(t)\|_X^2 \nonumber
\\
&
+C \sigma \|\nabla P_{\ne0}\omega^\nu(t)\|_{X}^2 \|\omega^\nu(t)\|_{L_2} +C\sigma \|\omega^\nu(t) \|_{L_2}\|\nabla\omega^0(t)\|_{X}^2.\label{n1n}
\eal

Therefore, it remains to provide a uniform estimate of $\|\omega^\nu(t)\|_{L_2}$ upper bounded by $\|\omega^\nu(0)\|_{L_2}$. To do so,  we note   that the operator $1+\Delta^{-1}$ is positive  on the space $(I-P_{\mathcal K})L_2(\T)$  due to $\alpha >1$,  and employ the $L_2$ estimate (\ref{L_2est1}) to obtain
\bal\min\{\frac{\alpha^2-1}{\alpha^2}, \frac34\} \| (I-P_{\mathcal K})\omega^\nu(t)\|_{L_2}^2&\le \|\omega^\nu(t)\|^2_X
\le \|\omega^\nu(0)\|^2_X
\le\|\omega^\nu(0)\|_{L_2} ^2.\label{nnn000}
\eal
Moreover, integrating by parts, we see  that
\be \lefteqn{\int_{\T}J(\Delta^{-1}P_{\mathcal K}\omega^\nu,\omega^\nu) P_{\mathcal K} \omega^\nu dxdy}
\\&=&
\int_{\T} \Big(\p_y(\omega^\nu\p_x\Delta^{-1}P_{\mathcal K}\omega^\nu)-\p_x(\omega^\nu\p_y\Delta^{-1}P_{\mathcal K}\omega^\nu)\Big) P_{\mathcal K}\omega^\nu dxdy=0.
\ee
Hence, taking the inner product of (\ref{PNS}) with $P_{\mathcal K}\omega^\nu$ and integrating by parts, we have
\bal \nonumber
\frac{d}{dt}& \|P_{\mathcal K}\omega^\nu(t)\|^2_{L_2}+2\nu \|P_{\mathcal K}\omega^\nu(t) \|^2_{L_2}
\\ \nonumber
 &= -2( J(\Delta^{-1}(I-P_{\mathcal K})\omega^\nu(t),(I-P_{\mathcal K})\omega^\nu), P_{\mathcal K}\omega^\nu(t))
\\ \nonumber
& \le 2\|\nabla \Delta^{-1} (I-P_{\mathcal K})\omega^\nu(t)\|_{L_2}\|(I-P_{\mathcal K})\omega^\nu(t)\|_{L_2}\|\nabla P_{\mathcal K}\omega^\nu(t)\|_{L_\infty}
\\ \nonumber
& \le C \|(I-P_{\mathcal K})\omega^\nu(t)\|_{L_2}^2\|P_{\mathcal K}\omega^\nu(t)\|_{L_2}
\\
&\le \nu \|P_{\mathcal K}\omega^\nu(t) \|^2_{L_2}+  C \frac1\nu\|\omega^\nu(0)\|_{L_2}^4.\label{nnn0000}
\eal
after the use of (\ref{nnn000}).
Integrating  (\ref{nnn0000}) and combining the resultant equation with  (\ref{nnn000}), we have
\bal\nonumber
\|\omega^\nu(t) \|^2_{L_2} &\le \Big(1+\frac{\alpha^2}{\alpha^2-1}+\frac43\Big)\| \omega^\nu(0)\|^2_{L_2}+  C \frac1\nu \int^t_0 e^{-\nu (t-s)} \| \omega^\nu(0)\|_{L_2}^4 ds
\\
&\le C\Big(\| \omega^\nu(0)\|^2_{L_2}+  \frac1{\nu^2} \| \omega^\nu(0)\|_{L_2}^4\Big)\le C \nu^2 d^2 \label{nnn11}
\eal
due the smallness assumption of  $\|\omega^\nu(0)\|_{L_2}$.

Therefore, by (\ref{nnn11}), the conservation of $\|\omega^0(t)\|_X^2$ from (\ref{L_2est1}) with $\nu=0$,  and the estimate \cite[Lemma 2.3]{Lin}
\be \|\nabla \omega^0(t)\|_X \le C(1+t)\|\nabla\omega^0(0)\|_X,
\ee
equation (\ref{n1n}) can be estimated as
\bal
\frac{d}{dt}\|\omega(t)\|^2_X\le& -\nu \|\nabla P_{\ne0}\omega^\nu(t)\|_X^2
+2\nu\|\nabla \omega^0(t)\|_X^2 +2\nu t^2\| \omega^0(t)\|_X^2 \nonumber
\\
&
+C \sigma \|\nabla P_{\ne0}\omega^\nu(t)\|_{X}^2 \|\omega^\nu(t)\|_{L_2} +C\sigma \|\omega^\nu(t) \|_{L_2}\|\nabla\omega^0(t)\|_{X}^2\nonumber
\\
\le &
 -\nu \|\nabla P_{\ne0}\omega^\nu(t)\|_X^2
+\nu C(1+t)^2\|\nabla \omega^0(0)\|_X^2 +2\nu t^2\| \omega^0(0)\|_X^2 \nonumber
\\
&
+C \sigma \nu d\|\nabla P_{\ne0}\omega^\nu(t)\|_{X}^2 +C\sigma \nu d  (1+t)^2\|\nabla\omega^0(0)\|_{X}^2
\nonumber
\\
\le &\nu C (1+t^2) \|\nabla \omega^0(0)\|_{X}^2
\label{n11n}
\eal
for $C$  independent of $\nu$ and for  $d$ small so that $C d \le 1$.
Hence, we obtain the validity of (\ref{LL}), after the integration of the previous equation.

The proof of Theorem \ref{th1} is complete.

\section{Enhanced and unenhanced dampings of exact solutions }

In this section, wee seek exact solutions, which damps in a  manner  not discussed in  Theorem  \ref{th1}.

Let $\omega$ solve   the Navier-Stokes equation  (\ref{NS}) with $\alpha >0$ and $a\in \{ 0,\,1\}$.
If $\omega$ satisfies   the  equation
\bbe \label{J}J(\Delta^{-1}\omega,\omega)=0,  \bee
which is valid when $\omega$ is a Laplacian  eigenfunction expressed as
\bbe \Delta \omega = \lambda \omega. \label{lam}\bee
Hence  (\ref{NS}) becomes the linear heat equation
\bbe\label{NSH} \p_t\omega =   \nu \Delta\omega  - \nu a\cos y,\bee
 and therefore $\omega$ is in the following form
\be \omega = -a\cos y +e^{\nu \Delta t}\Big( \omega(0) +a\cos y\Big) = -a\cos y+a\e^{-\nu t}\cos y +e^{\nu \Delta t}\omega(0). \ee

\subsection{Exact solutions of (\ref{NS}) with  $a=0,\, 1$}

For  any  function $\phi \in  P_0L_2({\T})$ or $P_{\ne0}\phi=0$, $\phi$ is a function independent of $x$. That is,
\be \phi=\sum_{n\ge 1}  (a_n \cos n y + b_n \sin n y)\in L_2(\T)
\ee
with constants $a_n$ and $b_n$.
 Therefore  (\ref{NS}) admits   the exact solution
\bal \label{s1}\omega &= -a\cos y+\e^{-\nu t} a\cos y +\e^{\nu \Delta t}\phi\nonumber
\\ &=-a\cos y+\e^{-\nu t} a\cos y + \sum_{n\ge 1}  e^{- n^2  \nu t}(a_n \cos n y + b_n \sin n y),
\eal
 since $\omega$ is independent of $x$ and hence  (\ref{J}) holds true.

 When $\alpha <1$, Kolmogorov flow losses stability, but the solutions (\ref{s1}) and (\ref{s2}) remain  in its  nonlinear stable manifold.

  The solution (\ref{s1})  is a parallel flow   moving in horizontal lines $y=$ constants, and is   in the kernel of the projection operator $P_{\ne0}$. Thus the solution (\ref{s1})  is  absent in the metastability estimate  in Theorem \ref{th1}, due to the energy estimate (\ref{nb1}). However, equation (\ref{s1}) shows that the stability can be enhanced if $\omega$  involves high frequency modes only.

\begin{figure}
 \centering
 \includegraphics[height=.35\textwidth, width=.49\textwidth]{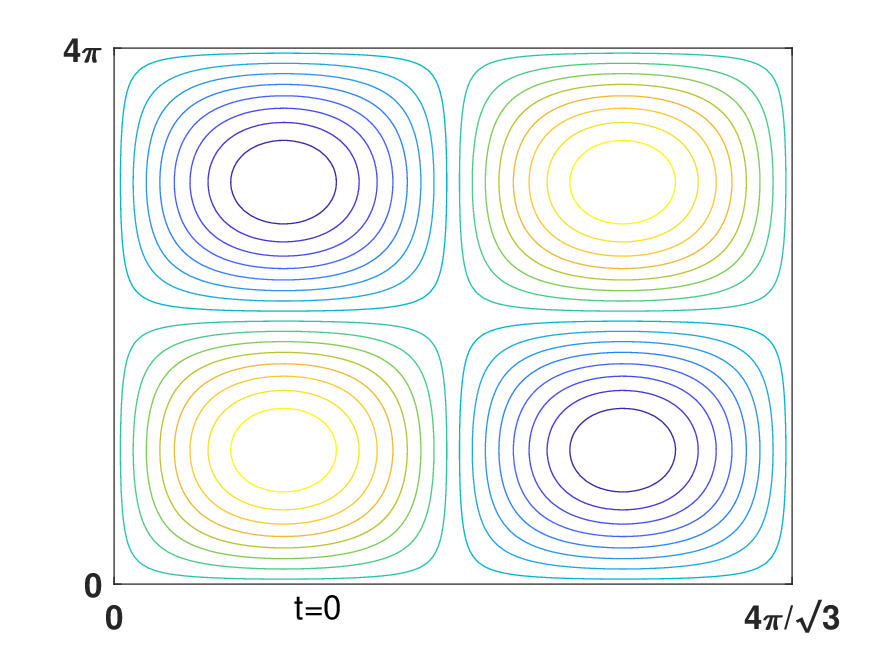}
 \includegraphics[height=.35\textwidth, width=.49\textwidth]{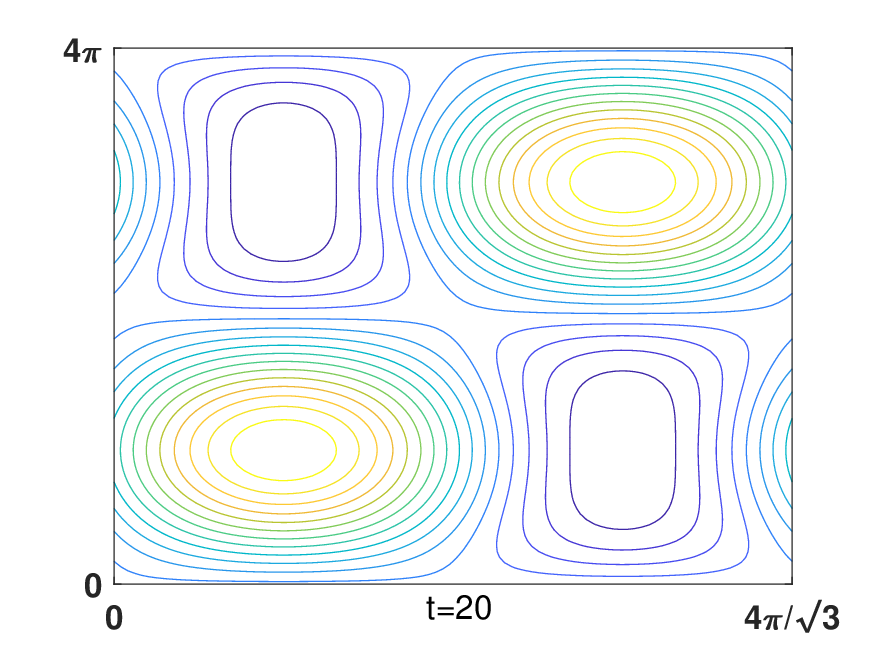}
 \includegraphics[height=.35\textwidth, width=.49\textwidth]{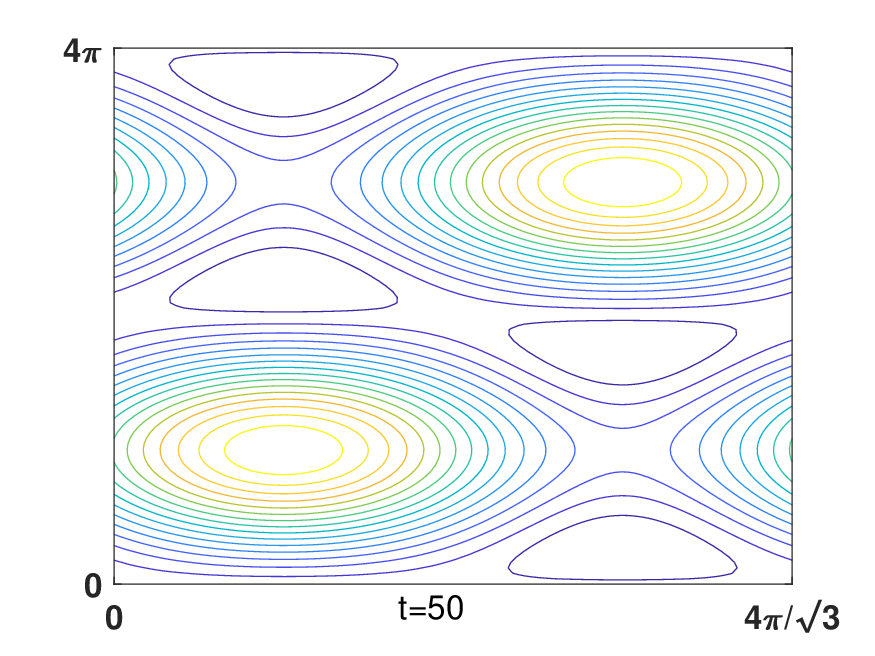}
 \includegraphics[height=.35\textwidth, width=.49\textwidth]{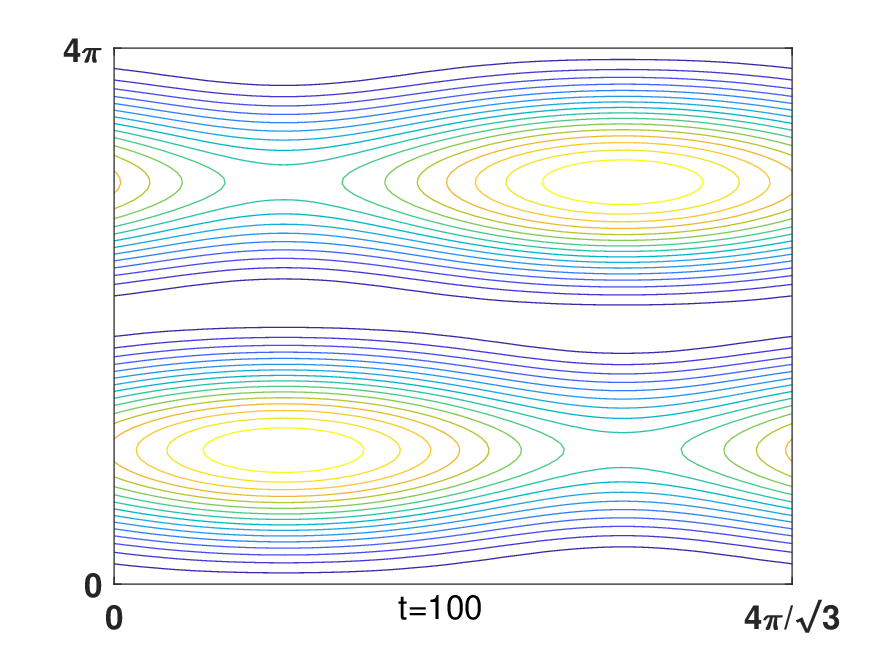}
\caption{Vorticity contours for  exact solution (\ref{s3}) when $\nu=0.01$.}
 \label{f1}
 \end{figure}

For nonparallel  flows, we consider (\ref{NS}) in the extended flat torus domain
\bbe \label{do}[0, 2\pi/\alpha)\times [0, 2\pi n) \mbox{ with } \alpha^2 + \frac1{n^2}=1 \mbox{ and } n\ge 1.\bee
We have the unenhanced damping  solution
\bal\nonumber
\omega=& -\cos y + e^{-\nu t}\cos y+ e^{-\nu t} \Big( c_1 \sin y +c_2 \cos y+c_3\sin \alpha x \sin \frac y{n} )
\\ &+ e^{-\nu t} \Big( c_4\cos \alpha x\cos \frac y{n} + c_5\sin \alpha x \cos \frac y{n} + c_6\cos \alpha x\sin \frac y{n} \Big)\label{s2}
\eal
 for any constants $c_i$, since (\ref{lam})  and hence  (\ref{J}) hold true.
For displaying purpose,  we choose a function from (\ref{s2}) as
\bbe
\omega= -\cos y +e^{-\nu t}\cos y + e^{-\nu t}\sin \frac{\sqrt{3}}2 x \sin \frac12 y    \label{s3}
\bee
in the torus domain $[0, \frac{4\pi}{\sqrt{3}})\times [0,4\pi)$. This solution  transforms initially from the four vortices state $\sin \frac{\sqrt{3}}2 x \sin \frac12 y$ to the horizontal parallel  flow $-\cos y $ (see Figure \ref{f1}  for  the transition). This transition  exhibits inverse cascade of two-dimensional fluid motion by enlarging two diagonal vortices, while other two vortices are compressed and then destroyed  into pieces.

\subsection{More exact solutions of (\ref{NS}) with  $a=0$}
For the case $a=0$, the basic stationary solution $-a\cos y$ becomes zero. We thus have more freedom to construct exact solutions.

To display further   solutions in $\T$ with $\alpha >0$,  we choose the following functions
\bal
\omega_1=& \sum_{k\ge 1}  e^{-\nu (n^2\alpha^2+m^2)k^2 t}(a_k\cos ( kn\alpha x\!+\!kmy) \!+\! b_k\sin ( kn\alpha x\!+\!kmy)) \in L_2(\T),\label{sol1}
\\
\omega_2 =& e^{-\nu (n^2\alpha^2+m^2) t} \Big(c_1\sin  n\alpha x \sin my+c_2\cos  n\alpha x \sin my\nonumber
\\
& \hspace{22mm}+c_3\sin  n\alpha x \cos my+c_4\cos  n\alpha x \cos my\Big), \label{sol3}
\\
\omega_3 =&  \omega_2+e^{-\nu j^2 t}( c_5\sin jy+ c_6\cos jy), \,\,\mbox{ when }  \alpha^2 n^2 +m^2 =j^2,\label{s10}
\\
\omega_4 =& \omega_2+  e^{-\nu i^2\alpha^2 t}( c_5\sin  i\alpha x +  c_6\cos  i\alpha x )\nonumber
\\ &+e^{-\nu j^2 t}( c_7\sin jy+ c_8\cos jy), \,\,\mbox{ when } \, \alpha^2n^2+m^2=\alpha^2 i^2=j^2,\label{s11}
\eal
 for any given positive  integers $n,\,m,\, i,\, j$  and  reals  $a_k,\, b_k,\, c_l$.

It is readily seen  that $\omega_2, \omega_3$ and $\omega_4$  solve (\ref{NS}) due to the validity of (\ref{lam}) and hence (\ref{J}).  To show  $\omega_1$ satisfying (\ref{J}) for $t\ge 0$, we use the notation
\be \phi_k=\cos ( kn\alpha x+kmy) \,\,\mbox{ and }\,\, \varphi_k=\sin ( kn\alpha x+kmy)
\ee
to obtain
\bal J&(\Delta^{-1}\omega_1,\omega_1)\nonumber
\\
&= -\sum_{k,k'\ge 1} e^{-\nu (k'^2+k^2)(\alpha^2n^2+m^2)t}\frac{kn\alpha (-a_k\varphi_k+b_k \phi_k) m k'(-a_{k'} \varphi_{k'} +b_{k'}\phi_{k'})}{n^2\alpha^2k^2+m^2k^2}\nonumber
\\
&+\sum_{k,k'\ge 1}e^{-\nu (k'^2+k^2)(\alpha^2n^2+m^2)t}\frac{km (-a_k\varphi_k+b_k \phi_k) n\alpha  k'(-a_{k'} \varphi_{k'} +b_{k'}\phi_{k'})}{n^2\alpha^2k^2+m^2k^2}=0.\nonumber
\eal
Therefore, $\omega_1$ solves (\ref{NS}) and is  a parallel  flow or a bar flow  moving along  the straight streamlines
\be n\alpha x +my = \mathrm{constants}.
\ee
For example, we take a function from (\ref{sol1}) as
\bal
\omega =& e^{-4 (\alpha^2 +1)\nu t} \sin (2 \alpha x+ 2y) +e^{-16 (\alpha^2 +1)\nu t} \cos (4 \alpha x+ 4y)\label{s4}
\eal
for $\alpha =\sqrt{6}/2$. 
\begin{figure}
 \centering
 \includegraphics[height=.35\textwidth, width=.48\textwidth]{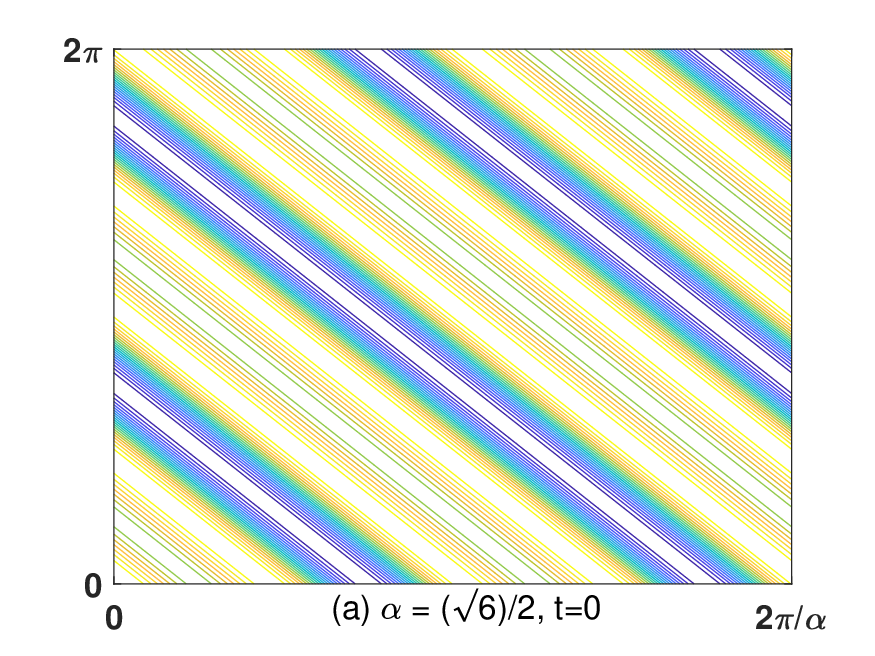}
 \includegraphics[height=.35\textwidth, width=.48\textwidth]{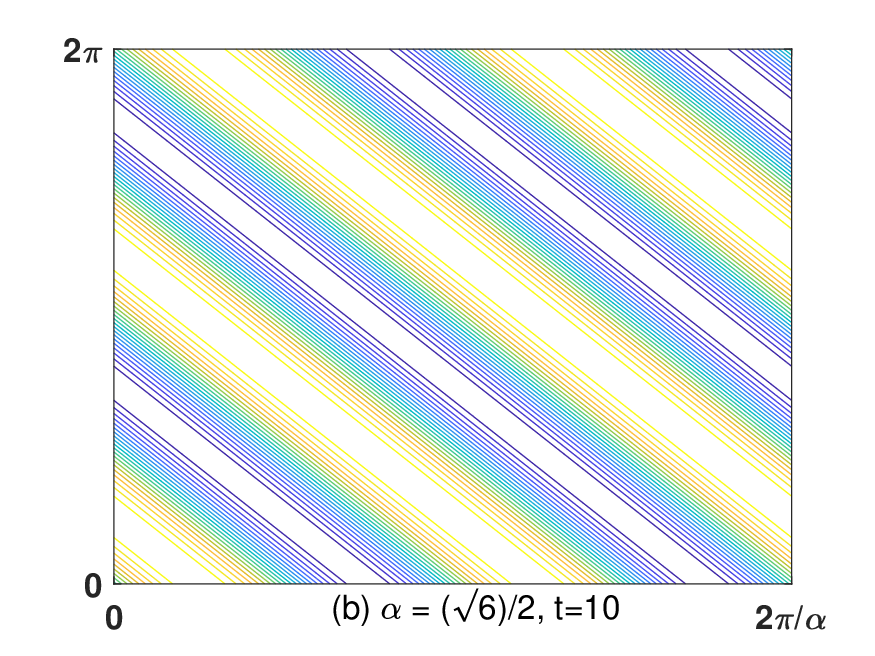}
 \includegraphics[height=.35\textwidth, width=.48\textwidth]{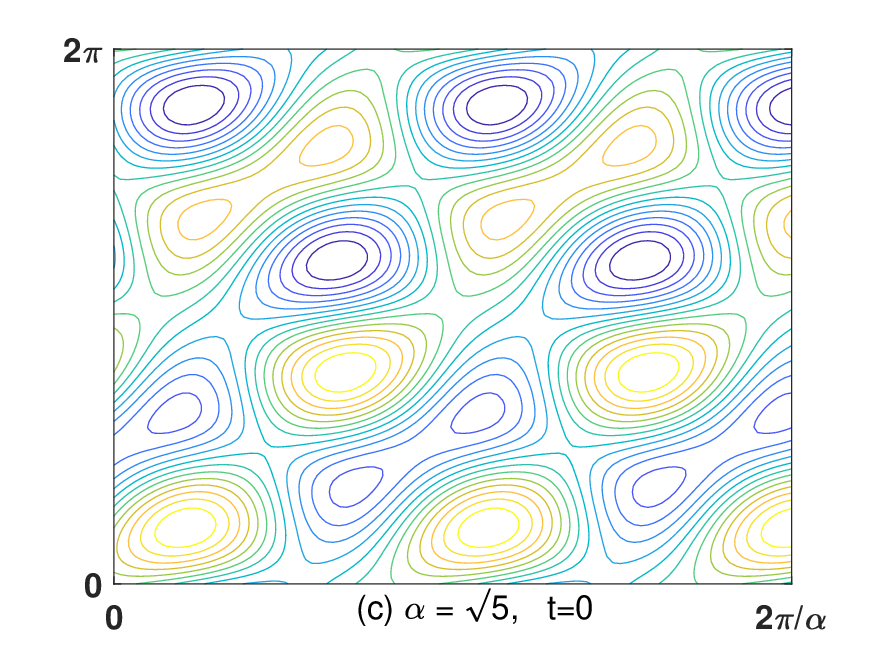}
 \includegraphics[height=.35\textwidth, width=.48\textwidth]{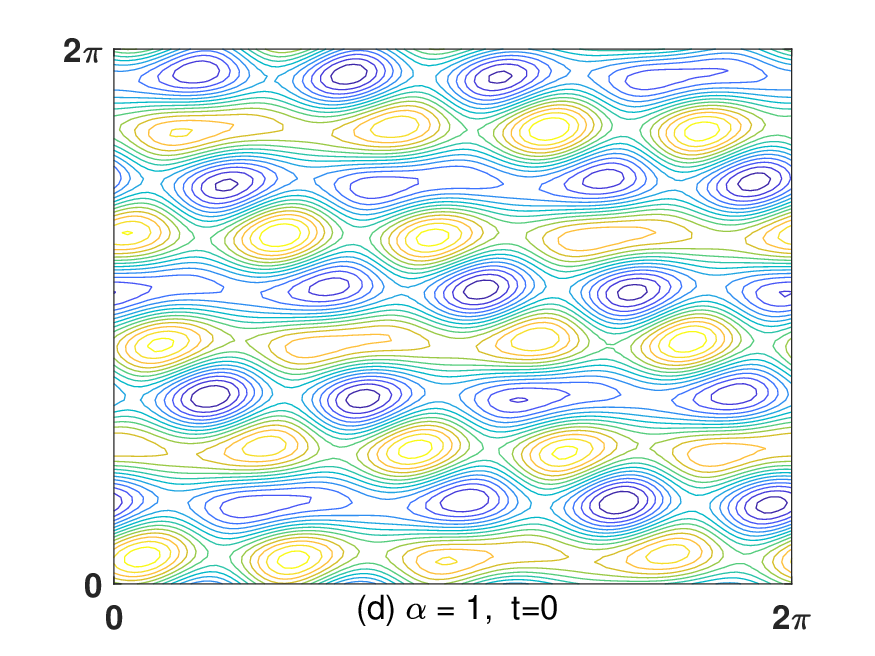}
\caption{ Vorticity contours of solutions (\ref{s4}) in (a)-(b), (\ref{s20}) in (c) and (\ref{s21}) in (d) for $\nu=0.01$.
}
 \label{f2}
 \end{figure}

For nonparallel flows, we choose a solution from (\ref{s10})   with $\alpha =\sqrt{5}$ expressed as
\bal
\omega&= e^{-(\alpha^2+4)\nu t} ( \sin (\alpha x)\sin (2y) +0.3 \cos(\alpha x)\cos(2y))+0.4\e^{-9\nu t}\sin (3y).\label{s20}
\eal
and adopt a solution from (\ref{s11})  with $\alpha=1$ in the following form
 \bal \omega =& e^{-(16 \alpha^2+9)\nu t} \Big(0.4 \sin (4\alpha x)\sin (3y) + 0.5 \cos (4\alpha x)\cos (3y))\nonumber
 \\ &+ \e^{-25 \nu t} \Big(\sin(5y)+0.3 \sin(5x)\Big).\label{s21}
 \eal

Vortex contours of enhanced damping  flows   (\ref{s4}), (\ref{s20}) and (\ref{s21}) are  displayed in Figure \ref{f2}.
  The  solutions (\ref{sol3})-(\ref{s11}) and their special forms  (\ref{s20})-(\ref{s21})  can be written as
 \bbe\label{fin}\omega(t)= e^{ -(a^2n^2+m^2)\nu t} \omega(0),\bee
which decays with $t$. The vorticity contours of $\omega(t)$  on the horizontal planes $z=$constants  are the same with those of $\omega(0)$ on the horizontal planes $z=e^{-(\alpha^2n^2+m^2) \nu t}$cosntants. Therefore, in Figure \ref{f2} (c)-(d), we only provide vorticity  contours at the initial state $t=0$, as the flow patterns  remain quasi-stationary  when $t$ grows.    On the other hand, when  $\nu$ is sufficiently small and $t$ is moderate, the exact solutions   evolve in a quasi-stationary manner and are close to their initial states, which solve the  stationary Euler equation.  Solutions (\ref{sol1})-(\ref{s11})   also exhibit  quasi-stationary states such as bar states (parallel  flows), dipole states and quadrupole states discussed by Yin {\it et al.} \cite{bar}. 

 It should be noted that the solution $e^{-(\alpha^2n^2 + m^2)\nu t} \sin (\alpha nx)\sin (my)$ is known as Taylor flow given by Taylor \cite{T}.

\begin{Remark}
If the fluid domain is the   horizontal channel  $\R \times [ 0, 2\pi]$  and the fluid motion is additionally driven by the up  moving boundary   in the following sense
\be \u|_{y=2\pi}  =(2\pi,0),\,\,\, \u|_{y=0}=(0,0), \ee
 the Navier-Stokes system (\ref{n1}) refers to  a  forced Couette flow problem and has  the exact parallel   solution
\be \u = \Big(y + a\sin y+ \sum_{n\ge 1}  e^{-\nu n^2 t} a_n\sin n y,\, 0\Big)\ee
for any coefficients $a_n$ so that  $\sum_{n\ge 1}n^2 a_n^2 <\infty$.

\end{Remark}

\

\noindent \textbf{Acknowledgement}. The present research is supported by the Shenzhen Natural Science Fund of China ( the Stable Support Plan Program No. 20220805175116001).

\

\end{document}